\input amstex
\input amsppt.sty
\magnification=\magstep1
\hsize=30truecc
\vsize=22.2truecm
\baselineskip=16truept
\TagsOnRight
\nologo
\pageno=1
\topmatter

\def\Z{\Bbb Z}

\def\C{\Bbb C}
\def\l{\left}
\def\r{\right}
\def\b{\bigg}
\def\bg{\bigg}
\def\({\b(}
\def\[{\b[}
\def\){\b)}
\def\]{\b]}
\def\colon{{:}\;}

\def\t{\text}
\def\f{\frac}
\def\mo{\roman{mod}}

\def\em{\emptyset}
\def\se {\subseteq}
\def\sp {\supseteq}
\def\sm{\setminus}

\def\eq{\equiv}
\def\cs{\cdots}
\def\ls{\leqslant}
\def\gs{\geqslant}
\def\al{\alpha}

\def\la{\lambda}
\def\Proof{\noindent{\it Proof}}
\def\Remark{\medskip\noindent{\it Remark}}
\def\Ack{\noindent {\bf Acknowledgment}}
\hbox {Math. Res. Lett. 11(2004), in press.}
\bigskip
\title Arithmetic Properties of Periodic Maps\endtitle
\author Zhi-Wei Sun\endauthor
\address Department of Mathematics (and Institute of Mathematical Science),
Nanjing University, Nanjing 210093,
The People's Republic of China\endaddress
\email {\tt zwsun}$\@${\tt nju.edu.cn}\  {\it Homepage}:\ {\tt http://pweb.nju.edu.cn/zwsun}\endemail
\abstract Let $\psi_1,\ldots,\psi_k$ be periodic maps from $\Z$
to a field of characteristic $p$ (where $p$ is zero or a prime). Assume that
positive integers $n_1,\ldots,n_k$ not divisible by $p$
 are their periods respectively. We show that $\psi_1+\cs+\psi_k$ is constant if
$\psi_1(x)+\cs+\psi_k(x)$ equals a constant for $|S|$ consecutive integers $x$
where $S=\bigcup_{s=1}^k\{r/n_s:\ r=0,\ldots,n_s-1\}$.
We also present some new results on finite systems of arithmetic
sequences.
\endabstract
\thanks 2000 {\it Mathematics Subject Classification}. Primary 11B25;
Secondary 11A07, 11A25, 11Y16, 68Q25.
\newline\indent
Supported by the Teaching and Research Award Fund for Outstanding
Young Teachers in Higher Education Institutions of MOE, and the
Key Program of the National Natural Science Foundation of P. R.
China.
\endthanks
\endtopmatter
\document
\heading {1. Introduction}\endheading
For $a\in\Z$ and $n\in\Z^+=\{1,2,3,\ldots\}$ we call
$$a(n)=a+n\Z=\{a+nx:\ x\in\Z\}$$
an {\it arithmetic sequence} with modulus $n$.
For a finite system
$$A=\{a_s(n_s)\}_{s=1}^k\tag1.1$$
of such sequences, the {\it covering function}
$w_A\colon\Z\to\Z$ given by
$$w_A(x)=|\{1\ls s\ls k:\ x\in a_s(n_s)\}|\tag1.2$$
is obviously periodic modulo the least common multiple
$[n_1,\ldots,n_k]$ of all the moduli $n_1,\ldots,n_k$. If
$w_A(x)\ls1$ for all $x\in\Z$ (i.e., $a_i(n_i)\cap a_j(n_j)=\em$
if $1\ls i<j\ls k$), then we say that $(1.1)$ is {\it disjoint}.
When $w_A(x)\gs1$ for all $x\in\Z$ (i.e.,
$\bigcup_{s=1}^ka_s(n_s)=\Z$), (1.1) is called
a {\it cover} of $\Z$.

A famous result of H. Davenport, L. Mirsky, D. Newman and R.
Rad\'o (cf. [NZ]) states that if (1.1) is a disjoint cover of $\Z$
with $1<n_1\ls\cs\ls n_{k-1}\ls n_k$ then we must have
$n_{k-1}=n_k$. In 1958 S. K. Stein [St] conjectured that if
$(1.1)$ is disjoint with $1<n_1<\cs<n_k$ then there exists an
integer $x\not\in\bigcup^k_{s=1}a_s(n_s)$ with $1\ls x\ls 2^k$. In
1965 P. Erd\H os [E2] offered a prize for a proof of his following
stronger conjecture (see [E1]): (1.1) forms a cover of $\Z$
 if it covers those integers from 1 to $2^k$.
 (The above $2^k$ is best possible because
$\{2^{s-1}(2^s)\}^k_{s=1}$ covers $1,\ldots,2^{k}-1$ but does not
cover any multiple of $2^k$.) In 1969--1970 R. B. Crittenden and
C. L. Vanden Eynden [CV1, CV2] supplied a long and awkward proof
of the Erd\H os conjecture for $k\gs20$.

 Let $m$ be a positive integer.
In [Su4, Su5] the author called (1.1) an {\it $m$-cover} of $\Z$
if $w_A(x)\gs m$ for all $x\in\Z$, and an {\it exact $m$-cover} of
$\Z$ if $w_A(x)=m$ for all $x\in\Z$.
Recently the author [Su10] found that $m$-covers of $\Z$
are closely related to subset sums in a field
and zero-sum problems on abelian groups.

Here is a result of [Su4, Su5] stronger than Erd\H os' conjecture:
{\it  $(1.1)$ forms an $m$-cover of $\Z$
if it covers
$|\{\{\sum_{s\in I}m_s/n_s\}:\, I\se\{1,\ldots,k\}\}|$
consecutive integers at least $m$ times, where
the given $m_1,\ldots,m_k\in\Z^+$ are
relatively prime to $n_1,\ldots,n_k$ respectively}.
(As usual the fractional part
of a real number $x$ is denoted by $\{x\}$.) In [Su5] the author
asked whether we have a similar result for exact $m$-covers of
$\Z$. The answer is actually negative, moreover there is no
constant $c(k,m)\in\Z^+$ such that  $(1.1)$ forms an exact
$m$-cover of $\Z$ whenever it covers $c(k,m)$ consecutive integers
exactly $m$ times. In fact, if (1.1) is an exact $m$-cover of $\Z$
then for any integer $N>1$ the system
$\{a_1(n_1),\ldots,a_k(n_k),0(N)\}$ covers $1,\ldots,N-1$ exactly
$m$ times but covers $0$ exactly $m+1$ times! (This observation is
due to the author's student H. Pan.)

For an assertion $P$ we adopt Iverson's notation
$$[\![P]\!]=\cases1&\t{if}\ P\ \t{holds},\\0&\t{otherwise}.
\endcases\tag1.3$$
Observe that $w_A(x)=\sum_{s=1}^k\psi_s(x)$ where
$\psi_s(x)=[\![n_s\mid x-a_s]\!]$ is periodic modulo $n_s$.

Our first result is completely new!

\proclaim{Theorem 1.1} Let $F$ be a field of characteristic $p$
where $p$ is zero or a prime. Let $n_1,\ldots, n_k$ be positive
integers not divisible by $p$, and let $\psi_1,\ldots,\psi_k$ be
maps from $\Z$ to $F$ with periods $n_1,\ldots,n_k$ respectively.
Then $\psi_1+\cs+\psi_k=0$ if $\psi_1(x)+\cs+\psi_k(x)=0$ for
$\sum_{d\in D}\varphi(d)$ consecutive integers $x$, where
$\varphi$ is Euler's totient function, $D=\bigcup_{s=1}^kD(n_s)$,
and $D(n)$ denotes the set of positive divisors of $n\in\Z^+$.
\endproclaim

\Remark\ 1.1. Clearly $\sum_{d\in D}\varphi(d)$ in Theorem 1.1
equals the cardinality of the set
$$\bigcup_{d\in D}\bg\{\f cd:\ 0\ls c<d\ \t{and}\ (c,d)=1\}\bg\}
=\bigcup_{s=1}^k\bg\{\f r{n_s}:\ r=0,1,\ldots,n_s-1\bg\},$$
where $(c,d)$ is the greatest common divisor of $c$ and $d$.
The result stated in the abstract is equivalent to Theorem 1.1 since
a constant can be viewed as a function on $\Z$ periodic mod 1.

\proclaim{Corollary 1.1} Let $w(x)$ be a function from $\Z$ to $\Z$ with period
$n_0\in\Z^+$. Then $w(x)$ is the covering function of $(1.1)$
if $w_A(x)=w(x)$ for
$|\bigcup_{s=0}^k\{0,1/n_s,\ldots,(n_s-1)/n_s\}|\ls
n_0+n_1+\cs+n_k-k$ consecutive integers $x$.
In particular, $(1.1)$ forms an exact $m$-cover of $\Z$
if it covers $|\bigcup_{s=1}^k\{r/n_s:\ r=0,\ldots,n_s-1\}|$
consecutive integers exactly $m$ times.
\endproclaim
\Proof. Let $D=\bigcup_{s=0}^kD(n_s)$. As
$$\psi(x):=w_A(x)-w(x)=-w(x)+\sum_{s=1}^k[\![n_s\mid x-a_s]\!]$$
vanishes for $|\bigcup_{s=0}^k\{r/n_s:\ r=0,\ldots,n_s-1\}|
=\sum_{d\in D}\varphi(d)$ consecutive integers
$x$, we have $\psi(x)=0$ for all $x\in\Z$ by Theorem 1.1.
When $n_0=1$ and $w(x)=m\in\Z^+$, this yields the latter result in Corollary 1.1.
\qed

\Remark\ 1.2. The problem whether a given $A=\{a_s(n_s)\}_{s=1}^k$
forms a cover of $\Z$ is known to be co-NP-complete. (See, e.g.
[GJ] and [T].) However, Corollary 1.1 indicates that we can check
whether system $A$ has a given covering function in
polynomial time! In 1997 the author [Su6] showed that if $(1.1)$
covers all the integers the same number of times then
$$\bg\{\sum_{s\in I}\f1{n_s}:\ I\se\{1,\ldots,k\}\bg\}
\sp\bigcup_{s=1}^k\l\{\f r{n_s}:\ r=0,\ldots,n_s-1\r\}.$$

\medskip
\noindent{\it Example 1.1}. Let (1.1) be an exact $m$-cover of
$\Z$, and let $n$ be an integer greater than $n_k$. Then the
system
$$A'=\{a_1(n_1),\ldots,a_{k-1}(n_{k-1}), a_k+n_k(n)\}$$
covers each of the consecutive integers $a_k+1,\ldots,a_k+2n_k-1$
exactly $m$ times but it does not cover $a_k$ or $a_k+2n_k$
exactly $m$ times. For example, $B=\{1(2),2(4),0(4)\}$ is a disjoint
cover of $\Z$, thus $B'=\{1(2),2(4),4(6)\}$ covers $1,\ldots,7$ exactly
once but it is not a disjoint cover. Note that the set
$\bigcup_{n\in\{2,4,6\}}\{r/n:\ r=0,\ldots,n-1\}$
just has 8 elements.

\proclaim{Corollary 1.2} Let $(1.1)$ be a system of arithmetic
sequences, and let $m$ be any integer  greater than
$k-f([n_1,\ldots,n_k])$. $($The function $f$ is given by $f(1)=0$
and $f(\prod_{i=1}^rp_i)=\sum_{i=1}^r(p_i-1)$ where
$p_1,\ldots,p_r$ are primes.$)$ Then there is an
$x\in\{0,1,\ldots,|S|-1\}$ such that $w_A(x)\not=m$ where
$S=\bigcup_{s=1}^k\{r/n_s:\ r=0,1,\ldots,n_s-1\}$.
\endproclaim

\Proof.  If (1.1) is an exact $m$-cover of $\Z$, then
$k\gs m+f([n_1,\ldots,n_k])$ by Corollary 4.5 of [Su7].
Thus, in view of the condition,
(1.1) does not form an exact $m$-cover of $\Z$ and
hence the desired result follows from Corollary 1.1. \qed

Our next theorem extends some earlier work in [Su4, Su5].

\proclaim{Theorem 1.2} Let $n_1,\ldots,n_k$ be positive integers,
and let $R_1,\ldots,R_k$ be finite subsets of $\Z$. For
$s=1,\ldots,k$, let $c_{st}$ lie in the complex field $\C$ for each
$t\in R_s$, and set
$$X_s=\bg\{x\in\Z:\  \sum_{t\in R_s}c_{st}e^{2\pi i\f
t{n_s}x}=0\bg\}.\tag1.4$$
If the system $\{X_s\}_{s=1}^k$ covers
$W$ consecutive integers at least $m$ times where $1\ls m\ls k$
and
$$W=\max\Sb I\se\{1,\ldots,k\}\\|I|=k-m+1\endSb
\bg|\bg\{\bg\{\sum_{s\in I}\f{r_s}{n_s}\bg\}:\  r_s\in R_s\bg\}\bg|
\ls\max\Sb I\se\{1,\ldots,k\}\\|I|=k-m+1\endSb\prod_{s\in I}|R_s|,\tag1.5$$
then it covers every integer at least $m$ times.
\endproclaim

\proclaim{Corollary 1.3} Let $(1.1)$ be a system of arithmetic
sequences, and let $m_1,\ldots,m_k$ be integers relatively prime
to $n_1,\ldots,n_k$ respectively. Let $l$ be any nonnegative integer with
$w_A(x)\gs l$ for all $x\in\Z$, and set
$$W_l=\max\Sb I\se\{1,\ldots,k\}\\|I|=k-l\endSb
\bg|\bg\{\bg\{\sum_{s\in J}\f{m_s}{n_s}\bg\}:\ J\se I\bg\}\bg|\ls 2^{k-l}.\tag1.6$$
Then the covering function $w_A(x)$ takes its minimum on every set
of $W_l$ consecutive integers.
\endproclaim
\Proof. Without loss of generality we may assume that $1\ls m_s\ls
n_s$ for all $s=1,\ldots,k$. As $m(A)=\min_{x\in\Z}w_A(x)\gs l$
and $W_l\gs W_{m(A)}$, it suffices to work with $l=m(A)$ below.

The case $l=k$ is trivial, so we let $l<k$. Set $c_{s0}=1$ and
$c_{sm_s}=-e^{-2\pi ia_sm_s/n_s}$ for $s=1,\ldots,k$. Since $m_s$
and $n_s$ are relatively prime,
$$X_s:=\l\{x\in\Z:\ c_{s0}e^{2\pi i\f 0{n_s}x}+c_{sm_s}e^{2\pi
i\f{m_s}{n_s}x}=0\r\}=a_s(n_s).$$
Applying Theorem 1.2 with $m=l+1$ and $R_s=\{0,m_s\}\ (1\ls s\ls
k)$, we immediately get the desired result. \qed

\Remark\ 1.3. (a) [Su9] contains some other interesting results on the
covering function of (1.1). (b) $W_l$ in (1.6) might be smaller than
its value in the case $m_1=\cdots=m_k=1$. Let $n_1=3,\ n_2=5$ and
$n_3=15$. Set
$$W_0(m_1,m_2,m_3)=\bg|\bg\{\bg\{\sum_{s\in J}\f{m_s}{n_s}\bg\}:
\ J\se\{1,2,3\}\bg\}\bg|$$
for $m_1,m_2,m_3\in\Z$. Then $W_0(1,1,2)=7<W_0(1,1,1)=8$.
\medskip

Our third theorem characterizes the least period of a covering
function.

 \proclaim{Theorem 1.3} Let $\la_s\in\C,\ a_s\in\Z$ and
$n_s\in\Z^+$ for $s=1,\ldots,k$. Then the smallest positive period
$n_0$ of the (weighted) covering function
$$w(x)=\sum_{s=1}^k\la_s[\![n_s\mid x-a_s]\!]$$
is the least $n\in\Z^+$ such that $\al n\in\Z$ for all
those $\al\in[0,1)$ with $\sum\Sb 1\ls s\ls k\\\al
n_s\in\Z\endSb\la_sn_s^{-1}e^{2\pi i\al a_s}\not=0$.
\endproclaim

\Remark\ 1.4. Under the condition of Theorem 1.3, it can be easily
checked that $\sum_{x=0}^{N-1}w(x)/N=\sum_{s=1}^k\la_s/n_s$ where
$N=[n_1,\ldots,n_k]$. If $w(x)=0$ for all $x\in\Z$, then $n_0=1$
and hence
$$\sum^k\Sb s=1\\\al n_s\in\Z\endSb\f{\la_s}{n_s}e^{2\pi i\al
a_s}=0\quad\ \t{for all}\ \al\in[0,1).\tag1.7$$ This was first
obtained by the author [Su2] in 1991 via an analytic method, and
the converse was proved in [Su3]. In [Su8] the author determined
those functions $f\colon\bigcup_{n\in\Z^+}\Z/n\Z\to\C$ such that
$\sum_{s=1}^k\la_sf(a_s+n_s\Z)$ only depends on the covering
function $w(x)$, this was announced by the author [Su1] in 1989.
\medskip

Let $l$ be a positive integer, and let
$$\Z^l=\{\vec x=\langle x_1,\ldots,x_l\rangle:\ x_1,\ldots,x_l\in\Z\}$$
be the direct sum of $l$ copies of the ring $\Z$. For $\vec x,\vec
y\in\Z^l$, we use $\vec x\mid\vec y$ to mean that $\vec y=\vec
q\vec x=\langle q_1x_1,\ldots,q_lx_l\rangle$ for some $\vec
q\in\Z^l$. A function $\Psi:\Z^l\to\C$ is said to be {\it periodic
modulo} $\vec n\in\Z^l$ if $\Psi(\vec x)=\Psi(\vec y)$ whenever
$\vec x-\vec y=\langle x_1-y_1,\ldots,x_l-y_l\rangle$ is divisible
by $\vec n$. For $x_1,\ldots,x_l\in\Z$, we also use $[x_t]_{1\ls
t\ls l}$ to denote the least common multiple of $x_1,\ldots,x_l$.

\proclaim{Theorem 1.4} Let $\la_s\in\C$, $\vec a_s\in\Z^l$ and
$\vec n_s\in(\Z^+)^l$ for $s=1,\ldots,k$ where $l\in\Z^+$. Suppose
that the function
$$w(\vec x)=\sum_{s=1}^k\la_s[\![\vec n_s\mid\vec x-\vec a_s]\!]\tag1.8$$
is periodic modulo $\vec n_0\in(\Z^+)^l$.  Let $\vec
d\in(\Z^+)^l$, $\vec d\nmid\vec n_0$ and
$$I(\vec d\!\ )=\{1\ls s\ls k:\ \vec d\mid \vec n_s\}\not=\em.$$
If $\sum_{s\in I(\vec d\!\ )}\la_s/(n_{s1}\cdots n_{sl})\not=0$, then
$$|I(\vec d\!\ )|\gs\bg|\bg\{\bg\{\sum_{t=1}^l\f{a_{st}}{d_t}\bg\}
\colon s\in I(\vec d\!\ )\bg\}\bg|\gs\min\Sb 0\ls s\ls k\\\vec
d\,\nmid\,\vec n_s\endSb\[\f{d_t}{(d_t,n_{st})}\]_{1\ls t\ls l}
\gs p(d_1\cdots d_l)$$ where we use $p(m)$ to denote the least
prime divisor of an integer $m>1$.
\endproclaim

\Remark\ 1.5. Theorem 1.4 is a generalization of the main result
of [Su2] which corresponds to the case $l=1$ and improves
the Zn\'am--Newman result [N].

\proclaim{Corollary 1.4} Let $\la_s\in\C\setminus\{0\}$, $\vec
a_s\in\Z^l$ and $\vec n_s\in(\Z^+)^l$ for $s=1,\ldots,k$ where
$l\in\Z^+$. Suppose that all those moduli $\vec n_s$ which are
maximal with respect to divisibility are distinct. Then the
function $w(\vec x)$ given by $(1.8)$ is periodic modulo
$\vec n_0\in(\Z^+)^l$ if and only if $\vec n_0$ is divisible by all the
moduli $\vec n_1,\ldots,\vec n_k$.
\endproclaim

\Proof. If $\vec n_s\mid\vec n_0$ for all $s=1,\ldots,k$, then the
function $w(\vec x)$ is obviously periodic mod $\vec n_0$.

Now suppose that $w(\vec x)$ is periodic modulo $\vec n_0$
but not all the moduli divide $\vec n_0$. Then there exists
a maximal modulus $\vec n_r$ with respect to divisibility
such that $\vec n_r\nmid\vec n_0$. By the condition,
$$I(\vec n_r):=\{1\ls s\ls k:\ \vec n_r\mid\vec n_s\}
=\{1\ls s\ls k:\ \vec n_s=\vec n_r\}=\{r\}.$$
On the other hand, by Theorem 1.4 we should have $|I(\vec n_r)|\gs
p(n_{r1}\cdots n_{rl})$. The contradiction ends our proof. \qed

\Remark\ 1.6.  Corollary 1.4 in the case $l=1$ was essentially
established by \v S. Porubsk\'y [P].

\heading{2. Proofs of Theorems 1.1--1.4}\endheading

\proclaim{Lemma 2.1} Let $c_1,\ldots,c_n$ lie in a field $F$, and
let $z_1,\ldots,z_n$ be distinct elements of $F\sm\{0\}$. If
$\sum_{j=1}^nc_jz_j^x$ vanishes for $n$ consecutive integers $x$,
then it vanishes for all $x\in\Z$.
\endproclaim
\Proof. Suppose that $\sum_{j=1}^nc_jz_j^{h+i-1}=0$ for every
$i=1,\ldots,n$ where $h\in\Z$. Since the Vandermonde determinant
$$\|z_j^{i-1}\|_{1\ls i,j\ls n}=\left|\matrix
1&1&\hdots &1\\
z_1&z_2&\hdots &z_{n}&\\
\vdots&\vdots&\ddots&\vdots\\
z_1^{n-1}&z_2^{n-1}&\hdots&z_{n}^{n-1}
\endmatrix\right|=\prod_{1\ls i<j\ls n}(z_j-z_i)$$
does not vanish, by Cramer's rule we have $c_jz_j^h=0$ and hence
$c_j=0$ for all $j=1,\ldots,n$. Therefore $\sum_{j=1}^nc_jz_j^x=0$
for any $x\in\Z$. \qed

\medskip
\noindent{\it Proof of Theorem 1.1}. As $p$ does not divide
$N=[n_1,\ldots,n_k]$, the equation $x^N-1=0$ has $N$ distinct roots in the
algebraic closure $E$ of the field $F$. The multiplicative group
$\{\zeta\in E:\ \zeta^N=1\}$ of order $N$ is cyclic, so $E$
contains an element $\zeta$ of multiplicative order $N$. For
$a\in\Z$ and $1\ls s\ls k$, we have the geometric series
$$\f1{n_s}\sum_{r=0}^{n_s-1}\zeta^{\f N{n_s}ar}=[\![n_s\mid a]\!].\tag2.1$$
Therefore
$$\align\sum_{s=1}^k\psi_s(x)
=&\sum_{s=1}^k\sum_{a=0}^{n_s-1} [\![n_s\mid a-x]\!]\psi_s(a)
\\=&\sum_{s=1}^k\sum_{a=0}^{n_s-1}\f1{n_s}\sum_{r=0}^{n_s-1}
\zeta^{\f N{n_s}(a-x)r}\psi_s(a)
\\=&\sum_{s=1}^k\f1{n_s}\sum_{a=0}^{n_s-1}\psi_s(a)
\sum\Sb 0\ls\al<1\\\al n_s\in\Z\endSb
\zeta^{\al N(a-x)}
\\=&\sum_{\al\in S}\l(\zeta^{-\al N}\r)^x\(\sum^k_{s=1}
\f{[\![\al n_s\in\Z]\!]}{n_s}
\sum_{a=0}^{n_s-1}\psi_s(a)\zeta^{\al Na}\),
\endalign$$
where $S$ is the set
$$\{\al\in [0,1):\ \al n_s\in\Z\ \t{for
some}\ 1\ls s\ls k\}=\bigcup_{s=1}^k\bg\{\f r{n_s}:\
r=0,\ldots,n_s-1\bg\}.$$ As those $\zeta^{-\al N}$ with $\al\in S$
are distinct, applying Lemma 2.1 we find that
$\sum^k_{s=1}\psi_s(x)=0$ for $|S|$ consecutive integers $x$ if
and only if $\sum_{s=1}^k\psi_s(x)=0$ for all $x\in\Z$. By Remark
1.1, $|S|=\sum_{d\in D}\varphi(d)$. This concludes the proof. \qed

\medskip
\noindent{\it Proof of Theorem 1.2}. Clearly an integer $x$ is
covered by $\{X_s\}_{s=1}^k$ at least $m$ times if and only if $x$
is covered by $\{X_s\}_{s\in I}$ for all $I\se\{1,\ldots,k\}$ with
$|I|=k-m+1$.

Now let $I\se\{1,\ldots,k\}$ and $|I|=k-m+1$. For any $x\in\Z$, we
have
$$\align \prod_{s\in I}\sum_{t\in R_s}c_{st}e^{2\pi i\f t{n_s}x}
=&\sum_{r_s\in R_s\ \t{for}\ s\in I}
\bg(\prod_{s\in I}c_{sr_s}\bg)e^{2\pi ix\sum_{s\in I}r_s/n_s}
\\=&\sum_{\theta\in R(I)}C_{I,\theta}e^{2\pi i\theta x}
\endalign$$
where $$R(I)=\bg\{\bg\{\sum_{s\in I}\f{r_s}{n_s}\bg\}:
\ r_s\in R_s\bg\} \ \ \t{and}\ \ C_{I,\theta}=\sum\Sb r_s\in R_s
\ \t{for}\ s\in I\\\{\sum_{s\in I}r_s/n_s\}=\theta\endSb
\prod_{s\in I}c_{sr_s}.$$
Since those $e^{2\pi i\theta}$ with $\theta\in R(I)$
are distinct, by Lemma 2.1 the system $\{X_s\}_{s\in I}$ covers
$|R(I)|$ consecutive integers $x$ if and only if it covers all
$x\in\Z$.

 In view of the above, we immediately obtain the desired result.
\qed

\medskip
\noindent{\it Proof of Theorem 1.3}. Let
$S=\{0\ls \al<1:\ \al n_s\in\Z\ \t{for some}\ 1\ls s\ls k\}$ and
$$T=\bg\{0\ls \al<1:\ c_{\al}
=\sum\Sb 1\ls s\ls k\\\al n_s\in\Z\endSb\f{\la_s}{n_s}
e^{2\pi i\al a_s}\not=0\bg\}.$$
For each $s=1,\ldots,k$ the arithmetical
function $\psi_s(x)=\la_s[\![n_s\mid x-a_s]\!]$ is periodic modulo
$n_s$. By the proof of Theorem 1.1, for any $x\in\Z$ we have
$$w(x)=\sum_{s=1}^k\la_s[\![n_s\mid x-a_s]\!]
=\sum_{\al\in S}e^{-2\pi i\al x}c_{\al}
=\sum_{\al\in T}e^{-2\pi i\al x}c_{\al}.$$

Let $n$ be the least positive integer such that $\al n\in\Z$ for
all $\al\in T$. By the above, $w(x)=w(x+n)$ for all $x\in\Z$. Thus
$n_0\mid n$.

If $T=\em$, then $n=1$ and hence $n_0=n$.
In the case $T\not=\em$, we have
$$0=w(x)-w(x+n_0)=\sum_{\al\in T}e^{-2\pi i\al x}(1-e^{-2\pi i\al n_0})c_{\al}$$
for every $x=0,\ldots,|T|-1$, and hence $(1-e^{-2\pi i\al
n_0})c_{\al}=0$ for any $\al\in T$ (Vandermonde). Now that $\al
n_0\in\Z$ (i.e., $e^{-2\pi i\al n_0}=1$) for all $\al\in T$, we
have $n_0\gs n$ and thus $n_0=n$.

The proof of Theorem 1.3 is now complete. \qed

\medskip
\noindent{\it Proof of Theorem 1.4}. Let $\vec c$ be any vector in $\Z^l$
with $\vec d\nmid\vec c\vec n_0$. Then, for some $1\ls r\ls l$ we have
$d_r\nmid c_rn_{0r}$
Note that $\vec n_0$ divides the vector $\langle0,\ldots,0,n_{0r},0,\ldots,0\rangle$.
For any $x_1,\ldots,x_{r-1},x_{r+1},\ldots,x_l\in\Z$, since
$$\sum_{s=1}^k\(\la_s\prod^l\Sb t=1\\t\not=r\endSb[\![n_{st}\mid x_t-a_{st}]\!]\)
[\![n_{sr}\mid x_r-a_{sr}]\!]=w(\vec x)$$
is periodic mod $n_{0r}$ as a function of $x_r$,
by Theorem 1.3 we must have
$$\sum^k\Sb s=1\\d_r\mid c_rn_{sr}\endSb
\(\la_s\prod^l\Sb t=1\\t\not=r\endSb[\![n_{st}\mid x_t-a_{st}]\!]\)
\f{e^{2\pi i (c_r/d_r)a_{sr}}}{n_{sr}}=0.$$
(Recall that $(c_r/d_r)n_{0r}\not\in\Z$.)

Let $J=\{1\ls s\ls k:\ d_r\mid c_rn_{sr}\}$
and $\la_s'=\la_sn_{sr}^{-1}e^{2\pi i a_{sr}c_r/d_r}$
for $s\in J$. Given $r'\in\{1,\ldots,l\}\sm\{r\}$ and $x_t\in\Z$ with $t\not=r,r'$,
we have
$$\align&\sum_{s\in J}\(\la_s'\prod^l\Sb t=1\\t\not=r,r'\endSb[\![n_{st}\mid x_t-a_{st}]\!]\)
[\![n_{sr'}\mid x_{r'}-a_{sr'}]\!]
\\&\qquad\quad=\sum_{s\in J}\la_s'\prod^l\Sb t=1\\t\not=r\endSb[\![n_{st}\mid x_t-a_{st}]\!]=0
\endalign$$
for all $x_{r'}\in\Z$. By applying Remark 1.4
$l-1$ times we finally obtain that
$$\sum^k\Sb s=1\\\vec d\,\mid \,\vec c\vec n_s\endSb\f{\la_s}{n_{s1}\cdots
n_{sl}}e^{2\pi i\sum_{t=1}^la_{st}c_t/d_t}=0.\tag2.2$$

Set $m=\min_{0\ls s\ls k,\, \vec d\,\nmid\,\vec
n_s}[d_t/(d_t,n_{st})]_{1\ls t\ls l}$. Clearly $m\gs p(d_1\cdots
d_l)$. Let $c$ be any positive integer less than $m$. For
$s=0,1\cs,k$ we have
$$\vec d\mid c\vec n_s\Leftrightarrow d_t\mid cn_{st}
\ \t{for all}\ t=1,\ldots,l
\Leftrightarrow\[\f{d_t}{(d_t,n_{st})}\]_{1\ls t\ls l} \ \bigg|\
c\Leftrightarrow\vec d\mid \vec n_s.$$ In other words, $\vec d\mid
c\vec n_s$ if and only if $s\in I(\vec d\,)$. (2.2) in the case
$\vec c=\langle c,\ldots,c\rangle$ yields that
$$\sum_{s\in I(\vec d\!\ )}\f{\la_s}{n_{s1}\cdots n_{sl}}
e^{2\pi ic\sum_{t=1}^la_{st}/d_t}=0.$$

Let $\Theta=\{\{\sum_{t=1}^la_{st}/d_t\}:\ s\in I(\vec d\,)\}$.
Suppose that $|\Theta|<m$. Then for each $c=1,\ldots,|\Theta|$ we
have
$$\align&\sum_{\theta\in\Theta}e^{2\pi ic\theta}\sum\Sb s\in I(\vec d\!\ )
\\\{\sum_{t=1}^la_{st}/d_t\}=\theta\endSb\f{\la_s}{n_{s1}\cdots n_{sl}}
\\=&\sum_{s\in I(\vec d\!\ )}\f{\la_s}{n_{s1}\cdots n_{sl}}
e^{2\pi ic\sum_{t=1}^la_{st}/d_t}=0.
\endalign$$
By Lemma 2.1 this holds for all integers $c$, in particular $c=0$:
$$\sum_{s\in I(\vec d\!\ )}\f{\la_s}{n_{s1}\cdots
n_{sl}}=0.$$
This directly contradicts one of the hypotheses,
whence $|\Theta|\gs m$. \qed

\bigskip

\Ack. The author is indebted to the referee for his valuable
suggestions to improve the presentation.

\enddocument